\magnification=1200
%
%

\def\H{{\cal H}_X}
\def\Ah{{\bf A}_h}
\def\A{{\bf A}}
\def\gb{{\cal G}}

\def\C{{\bf C}}
\def\P{{\bf P}}
\def\R{{\bf R}}
\def\P{{\bf P}}

\def\Z{{\bf Z}}

\def\u{{\cal U}}
\def\n{{\cal N}}
\def\p{{\cal P}}

\def\go{{{\cal G}_0}}
\def\g1{{{\cal G}_R}}

\def\G1{{{G}_R}}
\def\po{{{\cal P}_0}}
\def\lo{{{\cal L}_0}}

\def\ho{{{\cal H}_0}}
\def\hpo{{{\cal H}_{{\cal P}_0}}}

\def\lo{{{\cal L}_0}}
\def\hlo{{{\cal H}_{{\cal L}_0}}}

\def\l{{\cal L}}
\def\r{\rho}

\def\c{\gamma}

\def\lra{\longrightarrow}

\def\mapright#1{\smash{\mathop{\longrightarrow}\limits^{#1}}}

\font\mysmall=cmr8 at 8pt

\centerline{\bf ON THE DEFORMATION QUANTIZATION OF COADJOINT ORBITS}

\centerline{\bf OF  SEMISIMPLE GROUPS}
\bigskip
\centerline{R. FIORESI}
\bigskip
\centerline{Department of Mathematics, University of California}
\centerline{Los Angeles, CA 90095-1555, USA.}
\centerline{{\mysmall e-mail: rfioresi@math.ucla.edu}}
\bigskip
\centerline{M. A. LLED\'O
\footnote*{Supported  by EEC
under TMR contract ERBFMRX-CT96-0045, (Politecnico di \break Torino). } 
}
\bigskip
\centerline{Dipartimento di Fisica, Politecnico di Torino,}
\centerline{Corso Duca degli Abruzzi, 24, 10129 Torino and}
\centerline{Istituto Nazionale di Fisica Nucleare (INFN),}
\centerline{Sezione di Torino, Italy.}
\centerline{{\mysmall e-mail: lledo@athena.polito.it}}
\bigskip
\bigskip
\bigskip
{\centerline {TO THE MEMORY OF MOSHE FLATO}}

\bigskip
\bigskip
{\centerline {\bf Abstract}}
{\mysmall In this paper we consider the problem of deformation
quantization of the algebra of polynomial functions on 
coadjoint orbits of semisimple Lie groups. 
The deformation of an orbit 
is realized by taking the quotient of the universal enveloping
algebra of the Lie algebra of the given Lie group, by a suitable ideal.
A comparison with geometric quantization in the
case of SU(2) is done, where both methods agree.}
\bigskip

{\bf 1. Introduction.}

\bigskip

A system in classical mechanics is given by a symplectic manifold $X$ which we call
 phase space and a function on $X$, $H$, which we call Hamiltonian. The points in $X$
 represent possible states of the system, the commutative algebra 
$C^\infty(X)$ 
is the set of classical observables, corresponding to possible measurements on the system, and
 the integral curves of the hamiltonian vector field $X_H$ represent
the time evolution  of the classical system.

A quantization of the classical system $X$ has three ingredients [Be],

\noindent 1. A family of  noncommutative complex algebras $\Ah$ depending on a real
 parameter $h$, which we will identify with Planck's constant, satisfying
$$
\Ah\mapsto \A=C^\infty(X)^\C \quad \hbox{when}\quad h\mapsto 0,
$$
or a suitable subalgebra of $C^\infty(X)^\C$ determined by physical requirements,
 but enough to separate the points of $X$.  $C^\infty(X)^\C$ denotes
the complexification of 
$C^\infty(X)$.

\noindent 2. A  family of linear maps $Q_h:\A\mapsto \Ah$, 
 called the {\it quantization maps} satisfying
$$
{Q_h(F)*_h Q_h(G)-Q_h(G)*_h Q_h(F)\over h}\mapsto \{F,G\}\quad \hbox{when}\quad h\mapsto 0.
 $$
where $\{,\}$ is the Poisson bracket in $\A$ (extended by linearity).

\noindent 3. A representation of $\A_h$ on a Hilbert space $\H$,
$R:\A_h\mapsto \hbox{End}(\H)$.
 The real functions in $\A_h$ (belonging to $C^\infty(X)$) are mapped into hermitian operators.

 The elements of $\A_h$ are the quantum observables and the rays in $\H$ are the states of the
 quantum system. Not every possible realization of $\A_h$ on a Hilbert space $\H$ satisfies
 the physical requirements for the quantum system, since the set of rays of $\H$ should be in
 one to one correspondence with the quantum  physical states. So further requirements should be
 imposed on $\H$.

\medskip

A first step to find a quantization of a physical system is the
construction of a formal deformation of the Poisson
algebra  classical observables [BFFLS]. In general, formal deformations do not
present a closed solution to the quantization problem. One needs
to see if it is possible to specialize the deformation to an interval
of values of the formal parameter $h$ (including 0, so the limit
$h\mapsto 0$ is smooth), besides constructing the Hilbert
space where this algebra is represented. Nevertheless having a formal
deformation is a powerful technical tool in the process of quantization.

A first approach to this problem appears in [Be].
Berezin explicitly 
computes $*$-products for K\"ahler manifolds that are homogeneous
spaces. His approach provides an explicit integral formula for a $*$-product
where h is a real number. In [RCG] a geometric construction of Berezin's
quantization is performed.
 
Later  De Wilde and Lecomte [DL]  and Fedosov [Fe] separately, constructed
and classified formal $*$-products on generic symplectic  
manifolds. Etingof and Kazhdan [EK] proved the existence of a formal
deformation for another class of Poisson manifolds, the Poisson-Lie
groups. Finally,  Kontsevich  [Ko] proved the existence of an
essentially unique formal $*$-product on general Poisson manifolds.

More recently Reshetekhin and Taktajan [RT], starting from Berezin's
construction, were able to give an explicit integral formula for the
formal $*$-product on K\"ahler manifolds.

\medskip
It is our purpose to study the deformation quantization of coadjoint
orbits of semisimple Lie groups. In [ALM] it has been proven that a
covariant $*$-product exists on the orbits of the coadjoint orbit that
admit a polarization.  We will consider the algebra of polynomials on 
coadjoint orbits. 
In the above mentioned works $*$-products are given on $C^{\infty}$
functions, however there is no guarantee that there is a subalgebra of 
functions  that is closed under it. Instead, we will obtain both
a formal deformation and a deformation for any real value of $h$
for the subalgebra of polynomial functions.

In [Ko] Kontsevich briefly describes the algebra of polynomials
over the dual of the Lie algebra (a Poisson manifold) as a special case
of his general formula for $*$-product on Poisson manifolds (this
special case was known long before [Ho1] [Gu]). 
He does not however consider the restriction of those polynomials to a
coadjoint orbit submanifold and, as he points out later, the knowledge
of  $*$-product on a certain domain is far from giving knowledge of
$*$-product on subdomains of it. The formulation of a star product on
some coadjoint orbits using this deformation of the polynomial algebra 
was investigated 
in the series of works [CG] [ACG] and [Ho2] (and references inside). 

Our approach  starts also from the fact that the universal enveloping algebra of a complex
semisimple
Lie algebra is the deformation quantization of the polynomial algebra on
the dual Lie algebra. By quotienting by a suitable ideal we get a 
deformation quantization of the polynomial algebra on a regular coadjoint
orbit. Using some known facts on real and complex orbits this gives us
 a deformation quantization on the  regular orbits of
compact  semisimple Lie groups. No selection of ordering rule is needed
for the proof, which means that we obtain a whole class of star products
on the orbit. A proof of the analiticity of the deformation  in the deformation parameter
is provided here, and  the convergence of the deformed
product for polynomials on the orbit is obtained. More general cases, as regular orbits of
noncompact Lie groups, involve some subtleties that are partially
explored in Section 2. Further developments  will be given in a subsequent
paper. Also, the extension of the proof to non regular (although still
semisimple) orbits is non trivial.

Our construction has the advantage that it is given in a coordinate
independent way. Also the symmetries and its possible representations
are better studied in this framework.  
The formal deformation is realized using a true deformation
of the polynomials on the complex orbit. We obtain the deformation quantization as 
a non commutative algebra
depending on a formal parameter $h$ containing a subalgebra in which $h$
can be specialized to any real value. 

\medskip
Geometric quantization is another approach to the problem. The elements of the quantum system
 are constructed using the geometric elements of the classical system. 
(For an introduction to
geometric quantization, see for example  [Pu] and references inside). 
In the case when the
 phase space is $\R^{2n}$, a comparison between both procedures, deformation and geometric
 quantization has been established [GV]. Less trivial systems, as coadjoint orbits, have been
 the subject of geometric quantization. The guiding principle is
the preservation of the symmetries of the
 classical system after the quantization. The idea of finding a unitary
 representation of the symmetry group naturally attached to the coadjoint orbit is known as the
 Kirillov-Kostant orbit principle. The action of the group on the Hilbert space of
the representation should be induced by the action of the group on the phase space as
 symplectomorphisms. The algebra of classical observables should be substituted by a
 noncommutative algebra and the group should act also naturally by conjugation on this algebra.

The procedure we used in constructing the formal deformation, that is
assigning an ideal in the enveloping algebra to the coadjoint orbit,
makes the comparison with geometric quantization easier. In Section 4 we
 show that in the special case of SU(2) there is an
isomorphism between our deformation quantization and the algebra of
twisted differential operators that appears in geometric quantization.
\medskip
The organization of the paper is as follows. In Section 2  we make a review of the algebraic
 properties of the coadjoint orbits on which our method of deformation  is based.
In Section 3 we prove the existence  of the deformation and describe it explicitly in terms
of a quotient of the enveloping algebra by an ideal. In Section 4 we make a comparison of our
 results with the results of geometric quantization for a particularly simple case, the coadjoint
orbits of SU(2).

\bigskip
{\bf 2. Algebraic Structure of Coadjoint Orbits of Semisimple Lie Groups.}
\bigskip

Let $\G1$ be a real Lie group and $\g1$ its Lie algebra. The coadjoint
 action of $\G1$ on $\g1^*$ is given by
$$
<\hbox{Ad}^*(g)\lambda ,Y>=<\lambda,\hbox{Ad}(g^{-1})Y> \quad \forall\;
g\in \G1,\quad \lambda\in \g1^*,\quad Y\in \g1.
$$
We will denote by $C_{\G1}(\lambda)$ (or simply $C_\lambda$ if $\G1$ can 
be suppressed without confusion) the orbit of the point 
$\lambda\in \g1^*$ under the coadjoint action of $\G1$.

Consider now the algebra of $ C^\infty$ functions on $\g1^*$, $ C^\infty(\g1^*)$. We can turn
 it into a Poisson algebra with the so called Lie-Poisson structure
$$
\{f_1,f_2\}(\lambda)=<[(df_1)_{\lambda},(df_2)_{\lambda}],\lambda>,
\qquad f_1,f_2 \in C^{\infty}(\g1^*), \quad \lambda \in \g1^*.
$$
If $f \in  C^{\infty}(\g1^*)$, $(df)_{\lambda}$ is a map from $\g1^*$ to $\R$,
so it can be regarded as an element of $\g1$ and $[\;,\;]$ is the Lie bracket in
$\g1$. 
By writing the Poisson bracket in linear  coordinates, it is clear that
$\R[\g1^*]$, the ring of polynomials on $\g1^*$, is closed under the Poisson
 bracket.

The Hamiltonian vector fields define an integrable distribution on $\g1^*$ whose integral
 manifolds (the symplectic leaves) are precisely the orbits of the coadjoint action. So all the
 coadjoint orbits are symplectic manifolds with the symplectic structure inherited from the
 Poisson structure on $\g1^*$.

\medskip

 Let $G$ be a connected complex, semisimple Lie group and $\gb$ its Lie algebra. We wish to
 describe the coadjoint orbits of different real forms of $G$. 
We can identify
$\gb$ and $\gb^*$ by means of the Cartan-Killing form, so we will work with the
adjoint action instead. 
We denote by  $\G1$ an arbitrary real form of $G$, and $\g1$ its Lie algebra.

 We start with the adjoint orbits of the complex group $G$ itself. Let $Z_s\in \g1\subset
 \gb$ be a semisimple element. The orbit of $Z_s$ in $\gb$ under $G$ will be
 denoted by $C_G(Z_s)$. It is well known that this orbit is a smooth complex
algebraic variety defined over
 $\R$ [Bo]. That means that the real form of $C_G(Z_s)$,  $C_G(Z_s)(\R)=C_G(Z_s)\cap\g1$ is a
 real algebraic variety. If $\G1$ is compact,
$C_G(Z_s)(\R)$ coincides with the real orbit $C_{\G1}(Z_s)$.
In general
$C_G(Z_s)(\R)$ is the union of several real orbits
 $C_{\G1}(X_i), i\in I$ for some finite set of indices $I$ [Va2]. 
Hence the real orbits are not always algebraic varieties. 
We will give one of such examples later. Still, the algebraic structure of the closely
 related manifold $C_G(Z_s)(\R)$ will be useful for the quantization.

 The algebra that we want to deform is the polynomial ring on
the complex orbit. 
When $C_G(Z_s)(\R)$ consists of one real orbit, 
the complex polynomial ring is the
complexification of the polynomial ring on the real orbit.
In this case, giving a formal deformation defined over $\R$
of the polynomial ring of the complex orbit is completely equivalent to
give a formal deformation of the polynomial ring of the real orbit.

In general $I$
will have many elements. One can always consider the algebra of
polynomials on $C_G(Z_s)(\R)$ and restrict it to each of the connected
components. The $*$-product we obtain can also be defined on the algebra
of restricted polynomials without ambiguity, so we have a deformation of
certain algebra of functions on the real orbit. Interesting subalgebras
of the restricted polynomials that still separate the points of the
real orbit could be found, being also closed under the $*$-product. We will
 see such kind of construction in an example. 

\medskip

We summarize now the classification of real  coadjoint orbits  [Va2] [Vo]. 
The easiest situation is when
$\G1$ is a compact group. 
In this case the orbits are real algebraic varieties defined by the 
 polynomials on $\gb$, invariant 
with respect to the coadjoint action. 
These invariant polynomials (or Casimir
polynomials) are in one to one correspondence with polynomials on the Cartan subalgebra that are
 invariant under the Weyl group. So every point in a Weyl chamber
determines a value of the invariant polynomials, and hence, an adjoint orbit.

 The general case is a refinement of this particular one. We will consider only orbits that
 contain a semisimple element $Z_s\in \g1$. There are two special cases: the elliptic orbit,
when the minimal polynomial of the
element $Z_e$ has only purely imaginary eigenvalues, 
and the hyperbolic orbits, when 
the minimal polynomial of $Z_h$
has only real eigenvalues. The general case $Z_s=Z_h+Z_e$ can be understood in terms of the special
cases. 

 Let us denote by $U$ a  compact real form of $G$ and $\u$ its Lie algebra, while $G_0$ and $\go$
 denote a non compact form and its Lie algebra. 
The involution $\theta:\go\mapsto \go$ induces the Cartan decomposition $\go=\lo
+\po$, and $\u=\lo+i\po$. $K$ is a maximal compact subgroup of $G_0$ with Lie algebra
 $\lo$.
We denote by $\hpo$ the maximal abelian subalgebra of $\po$ and  by $\hlo$ a CSA of $\lo$.
 $W(G_0,\hlo)$ and $W(G_0, \hpo)$ will denote the  Weyl
groups corresponding to the root systems
of $K$ ($W(G_0,\hlo)$) and the restricted root system of $G_0$ ($W(G_0, \hpo)$).

The set of hyperbolic orbits is in one to one
correspondence with the 
 set of orbits of $W(G_0,\hpo)$ on $\hpo$, while the set of  elliptic 
 orbits is in one to one correspondence  with the set of orbits of $W(G_0,\hlo)$ on
$\hlo$. In summary, each point in the Weyl chamber of the corresponding
root system  determines a unique semisimple orbit and vice versa.
\medskip
{\bf Example 2.1}. {\it Orbits of SO(2,1)}.
\medskip
We want to show explicitly an example  where the real form of the complex orbit is the union
 of two real orbits. The value of the invariant polynomials in this case
doesn't completely determine a  real orbit.

 Consider the connected component containing the identity of the noncompact orthogonal group
  SO(2,1)= $\{3\times 3\; \hbox{real matrices} \; \Lambda / \Lambda^T\eta\Lambda=\eta\}$,
 where
$$
\eta=\pmatrix{1&0&0\cr
0&1&0\cr
0&0&-1} 
$$

The Lie algebra so(2,1) is given by $\hbox{so(2,1)= span}\{G,\tilde E,
\tilde F\}$, where
$$
G=\pmatrix{0&1&0\cr
-1&0&0\cr
0&0&0}
\quad
\tilde E=\pmatrix{0&0&0\cr
0&0&1\cr
0&1&0},
\quad
\tilde F=\pmatrix{0&0&1\cr
0&0&0\cr
1&0&0},
$$
with commutation relations
$$
[G,\tilde E]=\tilde F,\quad [G,\tilde F]=-\tilde E,\quad [\tilde
E,\tilde F]=-G.
$$
The involutive automorphism associated to this noncompact form of so(3) is
$\sigma(X)=\eta X\eta$ 
so the Cartan decomposition is given by
$\lo=\hbox{span}\{G\}$ and $\po=\hbox{span}\{\tilde E, \tilde F\}$. $\lo$ is the Lie algebra of
SO(2), the maximal compact subgroup, which in this case is abelian.

The only Casimir polynomial is given in the coordinates $X=x\tilde
E+y\tilde F+zG$ by $P(X)=x^2 +y^2 -z^2$.
The elliptic orbits are classified by the elements $\{tG, t\in\R-\{0\}\}$, so the equation
 describing this orbit is
$$
x^2+y^2-z^2=-t^2
$$

Notice that $t$ and $-t$ define the same equation (the same value for the Casimir), but
 they define different orbits. In fact, the solution of the equation above is a double sheeted
 hyperboloid, each of the sheets being a different orbit (inside the past and future cone
respectively).

Consider now the following automorphism of so(2,1) (in the ordered basis we gave before)
$$A=\pmatrix{-1&0&0\cr
0&-1&0\cr
0&0&1}.
$$
$A$ can in fact be written as $A=$Ad$(g)$ with $g$ an element in the
complexification of  SO(2,1). In fact,
$$
g=\pmatrix{-1&0&0\cr
0&1&0\cr
0&0&-1},
$$
belongs to SO(3), the compact real form. Acting on the CSA, span$\{G\}$,
it gives the only Weyl reflection  (the Weyl group of SO(3) is
$\{\hbox{Id},-\hbox{Id}\}$), so $g$ is a representative of the non
trivial element in the Weyl group of SO(3).

Notice that the CSA of the maximal compact subgroup SO(2) and of SO(3)
have the same dimension, but the automorphism $A$ is just the Weyl
reflection of SO(3) that is ``missing'' in SO(2). $A$ takes
a point in one sheet of the hyperboloid and sends it to the other
sheet, so $A$ is a diffeomorphism between the two real
orbits. 

 Consider now the subalgebra of polynomials on $\gb$ that are invariant under
 $A$ (since $A^2=\hbox{Id}$, $\{\hbox{Id},A\}$ is a subgroup of automorphisms of so(2,1)).
 It is easy to see that it is also a Poisson subalgebra. Moreover,
since the Casimir polynomial is invariant under $A$, it is also possible
to define a  subalgebra of the polynomial algebra  of the complex orbit.
It is defined over $\R$, since $A$ leaves the real form so(2,1)
invariant. This algebra is contained as subalgebra in the algebra of
polynomial functions over the real
orbit (by polynomial functions we mean polynomials in the ambient space
restricted to the orbit).

The implementation of such kind of procedure for more general cases is
still under study and  will be written elsewhere.

Hyperbolic orbits are classified by the Weyl chamber of the restricted root system. One can take
 $\hpo=\hbox{span}\{\tilde E\}$, then $\ho=\hbox{span}\{\tilde E\}$ so the only root is the
 restricted root. The
 Weyl chamber is $\{t\tilde E, t\in\R^+\}$, so the hyperbolic orbits are given by
$$
x^2 +y^2-z^2=t^2
$$
This is a single sheeted hyperboloid, so in this case the orbit is an algebraic manifold.

Finally we have the orbits in the light cone (nilpotent orbits) satisfying 
$$
x^2 +y^2-z^2=0
$$   
There are three of them, one for z=0, others for $z>0$ and  $z<0$, 
but we are not studying nilpotent orbits here.

\bigskip

{\bf 3. Deformation  of the polynomial algebra of
regular coadjoint orbits of semisimple groups.}

\bigskip
{\bf Definition 3.1}.  Given a real Poisson algebra $\P$, a {\it formal
 deformation} of $\P$
is an associative algebra $\P_h$
over $\R[h]$, where $h$ is a formal
parameter, with the following properties:

\noindent a. $\P_h$ is isomorphic to $\P[[h]]$ as a  $\R[[h]]$-module.

\noindent b. The multiplication $*_h$ in $\P_h$ reduces  mod($h$) to the one in $\P$.

\noindent c. $\tilde F *_h \tilde G -\tilde  G *_h \tilde F = h\{F,G\}$
mod $(h^2)$, where $\tilde F, \tilde G \in \P_h$ reduce to $F,G \in \P$ mod($h$)
and $\{\, ,\,\}$ is the Poisson bracket in $\P$.

If $X$ is a Poisson manifold and 
$\P=C^\infty(X)$ we call $\P_h$ a {\it formal deformation} of $X$. Some
authors also use the term {\it deformation quantization} of $X$.
\medskip
We can also speak of  the formal deformation of the
complexification $\A$ of a real Poisson algebra. The formal deformation of $\A$ will be an
associative algebra $\Ah$ with the same properties  (a),
(b) and  (c)  where $\R$ has been replaced by $\C$. We want to note here
that this doesn't convert the complexification of the symplectic
manifold $X$ in a real Poisson manifold of twice the dimension.
\medskip

We are going to describe first the  formal deformation of the polynomial
algebra on the complex orbit. 

In the first place we will consider $\C[h]$-modules, that is, we will
restrict the modules appearing on Definition 3.1 to be modules over $\C[h]$,
the algebra of the polynomials in the indeterminate $h$. 
This will
give us immediately the formal deformation by tensoring by
$\C[[h]]$. Notice that our formal deformation will contain a
subalgebra that can
be specialized to any value of $h\in \R$.

\medskip

Let $G$ be a complex semisimple Lie group of dimension $n$, $\gb$  
its Lie algebra and $U$ the enveloping algebra of $\gb$.
Let's denote by $T_A(V)$ the full tensor algebra of
a complex vector space $V$ over a $\C$-algebra $A$. 
Consider the  proper two sided ideal in $T_{\C[h]}(\gb)$
$$
{\cal L}_h=\sum_{X,Y \in \gb}
T_{\C[h]}(\gb) \otimes(X \otimes Y - Y \otimes X - h[X,Y])
\otimes T_{\C[h]}(\gb)
$$

We define $U_h=_{def} T_{\C[h]}(\gb)/\l_h$. $U_h$ can be interpreted in
the following way:

\noindent Let $\gb_h$ be the Lie algebra over $\C[h]$
$\gb_h=\C[h] \otimes_{\C} \gb$  with Lie bracket 
$$
[p(h) X,q(h) Y]_{h}=p(h)q(h)[X,Y]
$$
where $[\;,\;]$ and $[\;,\;]_h$ denote the brackets in $\gb$ and $\gb_h$
respectively. Then, $U_h$ is the universal enveloping algebra of
the algebra $\gb_h$.

We will denote with capital letters elements of the tensor algebras
and of $U_h$, while we will use lower case letters for
the elements of the polynomial algebra over $\gb^*$, $\C[\gb^*]$.
The product of two elements $A,B \in U_h$ will be written $AB$.
\medskip 

{\bf Proposition 3.2}. {\it  (Poincar\'e-Birkhoff-Witt theorem for $U_{h})$.
Let $\{X_1, \dots , X_n\}$ be a basis for $\gb$. Then
$$
1, X_{i_1}\cdots  X_{i_k} \qquad 1 \leq i_1 \leq \cdots \leq i_k \leq n 
$$
form a basis for $U_h$ as $\C[h]$-module.} 
\medskip

$U_{h}$ is a free $\C[h]$-module. In particular,
$U_h$ is torsion free. 
\medskip

{\bf Definition 3.3}.
Let $S(\gb)=T_\C(\gb)/{\l}$, with
$$
{\l}=\sum_{X,Y \in \gb}
T_{\C}(\gb) \otimes(X \otimes Y - Y \otimes X )
\otimes T_{\C}(\gb),
$$
be the symmetric algebra of $\gb$. The natural homomorphism
from $T_{\C}(\gb)$ to $S(\gb)$ is an isomorphism if restricted to 
the symmetric tensors. Let $\lambda$ be the inverse of such isomorphism.

The canonical isomorphism $\gb^{**} \cong \gb$, can be extended to an
algebra isomorphism 
$\C[\gb^{*}]\cong S[\gb]$ where $\C[\gb^{*}]$ denotes the polynomial
algebra over $\gb^*$. The composition of such isomorphism with $\lambda$
will be  called the {\it symmetrizer map}. 

Let $\{X_1,\dots , X_n\}$ be a basis for $\gb$ and
$\{x_1,\dots , x_n\}$ the corresponding basis for $\gb^{**} \subset \C[\gb^*]$. 
Then the symmetrizer map $\hbox{Sym} :\C[\gb^*]\longrightarrow  T_{\C}(\gb)$ is
given by
$$
\hbox{Sym}(x_1\cdots x_n)={1 \over p!}\sum_{s \in S_p}X_{s(1)}\otimes\cdots\otimes X_{s(p)}.
$$
where $S_p$ is the group of permutations of order $p$.
\medskip

Let $I\subset \C[\gb^*] $ be the set of  polynomials on $\gb^*$ invariant under 
the coadjoint action,
$$
I=\{p\in \C[\gb^*] \; | \;  p(\hbox{Ad}^*(g)\xi)=p(\xi) \quad 
\forall \xi \in \gb^*, \, g \in G\}.
$$
By Chevalley theorem we have that $I=\C[p_1,\dots , p_m]$, where
$p_1,\dots , p_m$ are algebraically independent homogeneous polynomials
and $m$ is the rank of $\gb$.

\medskip

{\bf Definition 3.4}.
We define a {\it Casimir element} in $T_{\C}(\gb)$ as the image
of an invariant polynomial under the symmetrizer map. 
Since $T(\gb) \subset T_{\C[h]}(\gb)$ 
Casimirs are also elements of $T_{\C[h]}(\gb)$.  
We call {\it Casimir element} in $U$ (respectively $U_h$) 
an element which is the image
of a {\it Casimir element in} $T(\gb)$ (respectively in $T_{\C[h]}(\gb)$)
under the natural projection.
\medskip
It is well known that the Casimir elements  lie in the center of $U$. We
want now to prove that they also lie in the center of $U_h$. 

Let's denote by $\tilde U_{h_0}$ the algebra
$U_h/((h-h_0){\bf 1})$, where $h_0 \in \C$, and by
$ev_{h_0}$ the natural projection $U_h \longrightarrow \tilde U_{h_0}$.

\medskip

{\bf Lemma 3.5}.
{\it
Let $P$ be a Casimir in $U_h$. Then $ev_{h_0}(P)$ is in the center of 
$\tilde U_{h_0}$.
}
\smallskip

{\sl Proof}. This is because $\tilde U_{h_0}$ is the
universal enveloping algebra of $\gb_{h_0}$, where $\gb_{h_0}$
is the complex Lie algebra coinciding with $\gb$ as vector space and
with bracket $[X,Y]_{h_0}=h_0[X,Y]$ where $[\;,\;]$ is the bracket in $\gb$.

\medskip

{\bf Theorem 3.6}. {\it The Casimir elements lie in the center of
$U_{h}$.}

\smallskip

{\sl Proof}. Let $P$ be a Casimir element and let $X_1,\dots , X_n$
be generators for $\gb$ hence for $\gb_h$. 
We need to show: $PX_i=X_iP$ for all $1 \leq i 
\leq n$. 
$$
PX_i-X_iP=\sum_{1\leq i_1 \leq \cdots \leq i_k \leq n}
u_{i_1\dots i_k}(h)X_{i_1}\cdots X_{i_k}
$$
Let us apply the $ev_{h_0}$ map.
$$
\eqalign{
&ev_{h_0}(PX_i-X_iP-\sum_{1\leq i_1 \leq \cdots \leq i_k \leq n} 
u_{i_1\dots i_k}(h)X_{i_1}\cdots X_{i_k})= \cr
&-\sum_{1\leq i_1 \leq \cdots \leq i_k \leq n} u_{i_1\dots
i_k}(h_0)X_{i_1}\cdots X_{i_k}=0
\cr}
$$
because by Lemma 3.5 $ev_0(PX_i-X_iP)=0$.
Since there are no relations among the standard monomials
$X_{i_1}\cdots X_{i_k}$ (Proposition 3.2) 
we have that $ u_{i_1...i_k}(h_0)=0$.
Since this is true for infinitely many $h_0$ and since $ u_{i_1\dots i_k}(h)$ is 
a polynomial we have that $u_{i_1\dots i_k}(h) \equiv 0$.

\medskip

We now restrict our attention to the regular coadjoint orbits, that
is the orbits of regular elements.
 We recall here the definition of a
regular element in $\gb^*$. Consider the characteristic polynomial of
$\hbox{ad}^*(\xi)$, $\xi \in \gb^*$,
$$
\hbox{det}(T\cdot {\bf 1}-\hbox{ad$^*$}(\xi))=\sum_{i \geq m} q_i(\xi)T^i.
$$
where $m=\hbox{rank}\gb^*$. The $q_i$'s are invariant polynomials. An element
$\xi\in \gb^*$ is regular if $q_m(\xi) \neq 0$. The regular elements are dense 
in $\gb^*$ and they are semisimple. 
In particular the regular elements in a Cartan subalgebra form the interior of
the Weyl chambers.

The orbits of  regular elements  are  orbits of maximal dimension $n-m$.
Observe also that 
the 0-eigenspace coincides with the centralizer of $\xi$, $Z_{\xi}$.
A semisimple element
$\xi$ is regular if and only if dim$(Z_{\xi})=m$.

\medskip

Let us  fix the coadjoint orbit $C_{\xi}$ of a regular element $\xi \in \gb^*$.
 The ideal of polynomials
vanishing on $C_{\xi}$ is given by
$$
I_0=(p_i-c_{i0},\; i=1,\dots, m), \qquad c_{i0} \in \C,
$$
where the $p_i$ have been defined above (see after Definition 3.3). 
$I_0$ is a prime ideal or equivalently the orbit $C_{\xi}$ 
is an irreducible algebraic variety.
 (In fact, the orbit of any semisimple element, regular or not, is an
irreducible algebraic variety [Ks]).

\medskip

Let's consider the Casimirs 
$\hat P_i=\hbox{Sym}(p_i)$,
where the $p_1,...,p_m$ are
generators for $I$ that satisfy Chevalley theorem. Let $P_i$
be the image of $\hat P_i$ in $U_h$. Define 
the two sided ideal generated by the relations $P_i-c_i(h),\; i=1,...,m$: 
$$
I_h= (P_i-c_i(h),\; i=1,...,m) \subset U_h 
$$
for $c_i(h)=\sum_j c_{ij}h^j$, $c_{ij} \in \C$ ($c_i(0)=c_{i0}$, the
constants appearing in the definition of $I_0$).

It is our goal to give a basis of the algebra $U_h/I_h$ 
as $\C{[h]}$-module. We need first a couple of lemmas.

\medskip

{\bf Lemma (3.7)}. {\it Let $\xi\in\gb^*$ be a  regular element of $\gb^*$
(or equivalently a
point in which the centralizer has dimension equal to the rank of $\gb^*$).
Then $(dp_1)_{\xi}$, ..., $(dp_m)_{\xi}$ are linearly independent.}

\smallskip

{\sl Proof}. See [Va3].
\medskip

{\bf Lemma (3.8)}. {\it
Let $r$  be a fixed positive integer and let all the notation be as above.
Let
$$
\sum_{1 \leq i_1\leq  \cdots \leq i_r \leq m} 
a_{i_1\dots i_r}(p_{i_1}-k_{i_1})\dots(p_{i_r}-k_{i_r})=0
$$
with $ a_{i_1\dots i_r}\in \C[\gb^*]$, $k_{i_1} \dots k_{i_r} \in \C$.
Then $a_{i_1\dots i_r} \in 
(p_{1}-k_{1},\dots,p_{m}-k_{m}) \subset \C[\gb^*]$.}

\smallskip

{\sl Proof}.
By lemma (3.7) we can choose local coordinates $(z_1, \dots , z_{n})$
in a neighborhood of $\xi$
so that $z_i=p_{i}-k_{i}$, $i=1, \dots , m$.  
Since $a_{i_1 \dots i_r}(z_1, \dots , z_{n})$
are analytic functions, we can represent them as power series in 
$z_1,\dots , z_n$:
$$
a_{i_1\dots i_r}(z_1,\dots , z_n)=\sum_{1\leq j_1\leq \cdots j_s\leq  n\atop
0\leq s} a_{i_1\dots i_r, j_1\dots j_s}
z_{j_1} \cdots z_{j_s}.
$$
This can be rewritten as:
$$
\eqalign{
a_{i_1\dots i_r}(z_1,\dots , z_n)=&
\sum_{m+1 \leq j_1 \leq\cdots j_s \leq n\atop 0\leq s} a_{i_1\dots i_r, j_1\dots j_s}
z_{j_1} \cdots z_{j_s}+
\cr 
 &\sum_{1 \leq l_1\leq \cdots l_t \leq n \atop
l_1<m,\,1\leq t}
a_{i_1\dots i_r, l_1\dots l_t}
z_{l_1} \cdots z_{l_t}\cr}.
$$
By substituting into the given equation we get:
$$
\eqalign{
&\sum_{1 \leq i_1\leq \cdots i_r \leq m}\sum_{m+1 \leq j_1 \leq \cdots j_s \leq n \atop 0\leq s}
a_{i_1\dots i_r, j_1\dots j_s}
z_{j_1} \cdots z_{j_s}z_{i_1} \cdots z_{i_r}+ \cr 
&\sum_{1 \leq i_1 \leq \cdots i_r \leq m}
\sum_{1 \leq l_1\leq \cdots l_t \leq n \atop
l_1<m ,1\leq t}
a_{i_1\dots i_r, l_1\dots l_t}
z_{l_1} \cdots z_{l_t}z_{i_1} \cdots z_{i_r}
=0\cr}
$$
Notice that, by the way the sums are defined, and being $r$ fixed,  both terms in the above
 equations have no monomials in common. This implies that
$$
\sum_{1 \leq i_1\leq \cdots i_r \leq m}\sum_{m+1 \leq j_1 \dots j_s \leq n}
a_{i_1\dots i_r, j_1\dots j_s}
z_{j_1} \dots z_{j_s}z_{i_1} \dots z_{i_r}=0
$$
from which
$$
a_{i_1\dots i_r, j_1\dots j_s}=0 \quad \forall \;1\leq i_1\dots i_r\leq m,
\quad m+1 \leq j_1\cdots j_s
$$
This implies
$$
a_{i_1\dots i_r}(z_1 \dots z_m) \in (z_1 \dots z_r).
$$
That is, locally
$$
a_{i_1\dots i_r}=\sum b_{i_1\dots i_rj}(p_j-k_j).
$$
So we have obtained that for all $\eta$ 
in a neighbourhood of $\xi$:
$$
a_{i_1\dots i_r}(\eta)-\sum b^j_{i_1\dots i_r}(\eta)(p_{j}(\eta)-k_{j})=0
$$
But since this function is algebraic and $C_{\xi}$ is irreducible 
this means that this function
is identically 0 on $C_{\xi}$. Hence the Lemma is proven.

\medskip
 
Let's consider the projection $\pi:U_h \lra
U_h/(h{\bf 1})
\cong S(\gb) \cong \C[\gb^*]$. 
We have that  $\pi(A)=\pi(B)$ if and only if $A \equiv B$ mod$h$.
To simplify the notation we will 
denote the element of $\C[\gb^*]$ corresponding to $\pi(A)$
by $a$ (same letter, but lower case), as we did for the Casimirs $P_i$ before. 

\medskip

{\bf Lemma 3.9}. 
{\it Let $k$ be a fixed integer and let
$$
\sum_{i_1\leq\cdots i_k\leq m} A_{i_1\dots
i_k}(P_{i_1}-c_{i_1}(h))\cdots (P_{i_k}-c_{i_k}(h))
\equiv 0 \qquad \hbox{mod} h
$$
where $A_{i_1\dots i_k} \in U_h$ and the $P_{i}$'s and $c_{i}(h)$'s
have been defined above.
Then
$$
\matrix{
\sum_{i_1\leq\cdots i_k\leq m} A_{i_1\dots
i_k}(P_{i_1}-c_{i_1}(h))\cdots (P_{i_k}-c_{i_k}(h))
=h\sum_{i_1\leq\cdots i_k\leq m} B_{j_1 \dots j_l,i_1\dots
i_k} 
\cr \cr
(P_{j_1}-c_{j_1}(h))\cdots (P_{j_l}-c_{j_l}(h))
(P_{i_1}-c_{i_1}(h))\cdots (P_{i_k}-c_{i_k}(h))
}
$$
}
\smallskip

{\sl Proof}. By induction on $N=$max$_{i_1\dots i_k}$deg$a_{i_1\dots i_k}$, where,
using the the convention above,
$a_{i_1\dots i_k}=\pi(A_{i_1\dots i_k})$.
Let $N=0$. We have:
$$
\sum a_{i_1\dots i_k}(p_{i_1}-c_{i_10})\cdots (p_{i_k}-c_{i_k0})=0
$$
with $a_{i_1\dots i_k} \in \C$. 
By Lemma (3.8) $a_{i_1\dots i_k} \in I_0$ hence $a_{i_1\dots i_k}=0$.
This implies that $A_{i_1\dots i_k}=hB_{i_1\dots i_k}$.

Let's now consider a generic $N$,
$$
\sum a_{i_1\dots i_k}(p_{i_1}-c_{i_10})\cdots (p_{i_k}-c_{i_k0})=0
$$
By Lemma (3.8)
$$
a_{i_1\dots i_k}=\sum_j a_{i_1\dots i_kj}(p_j-c_{j0})
$$
with max${}_{i_1\dots i_k}\hbox{deg} a_{i_1\dots i_kj} < N$. Again we have that
$$
A_{i_1\dots i_k}=\sum_j A_{i_1\dots i_kj}(P_j-c_{j}(h))+hC_{i_1\dots i_k}.
$$
Let's substitute $A_{i_1\dots i_k}$
$$
\sum A_{i_1\dots i_kj}(P_j-c_{j}(h))(P_{i_1}-c_{i_1}(h))\cdots (P_{i_k}-c_{i_k}(h))
\equiv 0 \qquad \hbox{mod} h.
$$
By induction we have our result.

\medskip

{\bf Lemma 3.10}.
{\it If $hF \in I_h$ then $F \in I_h$.}

\smallskip

{\sl Proof}.
Since $hF \in I_h$ and since the $P_i$ are central elements:
$$
hF = \sum A_i (P_i-c_i(h))
$$
We have $\sum A_i (P_i-c_i(h)) \equiv 0$ mod$h$. Hence, by Lemma 3.9 and
also by the fact that $U_h$ is torsion free
we have our result.

\medskip

We have shown that $U_h/I_h$ is a $\C[h]$-module without torsion. We are ready
now to show that it is a free module by explicitly constructing a
basis. 
Let's fix a basis $\{X_1,\dots, X_n\}$ of $\gb$ and let $x_1,\dots , x_n$ 
be the corresponding elements in $\C[\gb^*]$. With this choice $\C[\gb^*]
\cong \C[x_1,\dots , x_n]$.
Let $\{x_{i_1},
\dots , x_{i_k}\}_{(i_1,\dots , i_k) \in \A}$ be  a basis in
of $\C[\gb^*]/I_0$ as $\C$-module, where $\A$
is a set of multiindices
appropriate to describe the basis. In particular, we can take them such 
that $i_1 \leq \cdots\leq i_k$.

\medskip

{\bf Proposition 3.11}.
{\it The monomials $\{X_{i_1}\cdots X_{i_k}\}_{(i_1, \dots , i_k) \in \A}$
are linearly independent in  $U_h/I_h$.}

\smallskip

{\sl Proof}.
Suppose  that there exists  a linear relation among the \break $X_{i_1},\cdots  X_{i_k}$'s,
 $(i_1,\dots, i_k) \in \A$ and let  $G \in I_h$  be such relation, 
$$
G=G_0+G_1h+ \cdots , 
\qquad G_i \in \hbox{span}_{\C}\{X_{i_1}\cdots X_{i_k}\}_{(i_1 \dots  i_k)\in \A}.
$$
Assume $G_i=0$,  $i<k$, $G_k \neq 0$.
We can write $G=h^kF$, with
$$
F=F_0+F_1h+ \cdots , \qquad F_0 \neq 0
$$
Since $h^kF \in I_h$ by hypothesis, using  Lemma (3.10) we have that
 $F \in I_h$, that is 
$$
F = \sum A_i(P_i-c_i(h)),
$$ 
and reducing  mod $h$,
$$
f = \sum a_i(p_i-c_{i0}).
$$
This would mean that $f$ represents a non trivial
relation among the monomials \break 
$\{x_{i_1}\cdots x_{i_k}\}_{(i_1 ... i_k) \in \A}$
in $C[\gb^*]/I_0$, which is a contradiction, so the linear independence
is proven.

\medskip

We want to give a procedure to construct a basis on $\C[\gb^*]/I_0$  starting from
a set of generators of $\C[\gb^*]$,  $S=\{x_{i_1}\cdots x_{i_k}\}$ $\forall \; 1 \leq i_1 \leq
\cdots i_k \leq n$.
 As a linear space
$I_0=$span$_{\C}\{x_{i_1}\cdots x_{i_k}(p_i-c_i)\}$.
Every element of the set  that spans $I_0$ 
will provide one relation that will allow us to eliminate  at most 
one element of the set $S$.
We can choose to eliminate successively the greatest element with respect to lexicographic
ordering. This means that any monomial in $S$ will be expressed in
terms of monomials of degree less or equal to its degree.
\medskip

{\bf Remarks (3.12)} We want to make two remarks that will be used
later.

\smallskip

1. An arbitrary monomial
$x_{j_1}\cdots x_{j_r}$ in $\C[\gb^*]$ can be written as:
$$
x_{j_1}\cdots x_{j_r}=\sum_{k \leq r\atop (m_1,\dots ,m_k)\in\A} 
a_{m_1\dots m_k}^{j_1\dots j_r}x_{m_1}\cdots x_{m_k}
+\sum_{i,d_i+g_i \leq r}
b_i(p_i-c_i)
$$
where $b_{i}$ is polynomial of degree $g_i$, $d_i$=deg$p_i$ and $a_{m_1\dots m_k}^{j_1\dots
j_r}\in \C$.

\smallskip

2. Let $A \in U_h$, $A \neq 0$, 
$A \in \hbox{span}_{\C}\{X_{j_1}\cdots X_{j_p}\}_{p \leq r}$, 
$j_1\dots j_p$ not necessarily
ordered. If $A \equiv 0$ mod$h$, then $A=hB$,
$B \in \hbox{span}_{\C}\{X_{i_1}\cdots X_{i_p}\}_{p < r\atop i_1 \leq\dots \leq i_p}$.

\medskip
 Next proposition will show the generation, so we will have a basis. 
\medskip

{\bf Proposition 3.13}. {\it
The standard  monomials $\{X_{i_1}\cdots  X_{i_k}\}$ with ${(i_1,\dots , i_k) \in \A}$
generate $U_h/I_h$ as $\C[h]$-module.}

\smallskip

{\sl Proof}. 
By Proposition 3.2 (PBW theorem in $U_h$) it is sufficient to prove that

$$
X_{j_1}\cdots X_{j_r} \in \hbox{span}_{\C[h]}\{X_{i_1} \cdots X_{i_k}
\}_{(i_1,\dots , i_k) \in \A}
$$
where $1 \leq j_1 \leq \cdots j_r \leq n$ and $X_{j_1}\cdots X_{j_r}$
denotes also the projection onto $U_h/I_h$ of the standard monomial

We proceed by induction on $r$. For $r=0$ it is  clear. For generic $r$
we write (see Remark 3.12)
$$
x_{j_1}\cdots x_{j_r}=\sum_{k \leq r\atop (m_1,\dots , m_k) \in \A} 
a_{m_1\dots m_k}^{j_1\dots j_r}x_{m_1}\cdots x_{m_k}
+\sum_{i, d_i+g_i\leq r}
b_{i}(p_i-c_i)  
$$
Lifting this equation from the symmetric algebra to the enveloping
algebra we have 
$$
X_{j_1}\cdots X_{j_r}-\sum_{k \leq r\atop (m_1,\dots , m_k) \in \A} a_{m_1\dots m_k}^{j_1\dots
j_r}X_{m_1}\cdots X_{m_k}
-\sum_{i}B_i(P_i-c_i(h))=hB 
$$
where, by the remark 2 in 3.12,  $B \in \hbox{span}\{X_{i_1}\cdots X_{i_p}\}_{p < r}$.
Applying the induction hypothesis, we have our result.
\medskip

Let $\C_h[\gb^*]=\C[h] \otimes \C[\gb^*]$, $I_0'=\C[h] \otimes I_0$.
We are now ready to prove the following theorem:

\medskip
 
{\bf Theorem (3.14)}. {\it Let the notation be as above. We have that
$U_h/I_h$ has the following properties:

\smallskip

\noindent 1. $U_h/I_h$ is isomorphic to $\C_h[\gb^*]/I_0'$ as
a $\C[h]$-module.

\smallskip

\noindent 2. The multiplication in $U_h/I_h$ reduces $mod (h)$ to the one in
$\C[\gb^*]/I_0'$.

\smallskip

\noindent 3. If $FG-GF=hP$, $F, G, P \in U_h/I_h$,  then $p=\{f,g\}$, where $\{,\}$ is the Poisson
bracket on the orbit defined by $I_0$.} (We are using the same
convention, $f=\pi(F)$).

\smallskip

{\sl Proof.}

\noindent 1. It is a consequence of Propositions 3.11, 3.13.

\noindent 2. It is is trivial.

\noindent 3. This property is satisfied by  the multiplication
in $U_h$ and the Poisson bracket in $\C[\gb^*]$ (see [Ko], [CP], [Ki]).
The Poisson bracket in the $\C[\gb^*]/I_0$ is induced from the one in
$\C[\gb^*]$, it is enough to see that $p$ will not depend on the
representative chosen in $U_h/I_h$, which is trivial.
  
\medskip

It is now immediate to obtain the properties of Definition 3.1 when we
consider the extension of $\C[h]$ to $\C[[h]]$. We define 
$$
\matrix{
\C_{[h]}[\gb^*]=\C[[h]] \otimes \C[\gb^*] & I_{[0]}=\C[[h]] \otimes I_0
\cr \cr
U_{[h]}=\C[[h]] \otimes U_h & I_{[h]}=\C[[h]] \otimes I_h.
}
$$

\medskip

{\bf Theorem 3.15}  
{\it
$U_{[h]}/I_{[h]}$ is a formal deformation (or a deformation quantization) of
$\C_{[h]}[\gb^*]/I_{[0]}$.}

\medskip

We want to note here that whatever is the real form chosen, the deformed
algebra is defined over $\R$, provided $c_{ij} \in \R$. 
Care should be taken, nevertheless, in
choosing the appropriate generators of $I_0$ with real coefficients
and this is always possible ([Bo]). 

\medskip 

Finally we want to come back to Example 2.1 and exhibit the deformed algebra.

\medskip

{\bf Example 3.16}. Let $G=SL_2(\C)$. 
The standard basis for $\gb=sl_2(\C)$ is  $\{H,X,Y\}$ with commutation relations
$$
\matrix{
[H,X]=2X & [H,Y]=-2Y & [X,Y]=H.
}
$$
We identify $\gb$ and $\gb^*$ via the Cartan Killing form. The  only
independent   invariant polynomial is: 
$$
p={1\over 4}h^2 + xy
$$
or, in terms of the compact  generators
$$\eqalign{&E={1\over 2}(X-Y) \quad F={i/2}(X+Y) \quad G={i/2}H\cr
&p=-(e^2+f^2+g^2)\cr}$$

The orbit $C_\xi$ of the regular semisimple
element $\xi=\pmatrix{ia/2 & 0 \cr 0 & -ia/2}$ (see the fundamental
representation in the next section) has coordinate
ring $\C[h,x,y]/(e^2+f^2+g^2-a^2)$. So we have that
$$
U_{[[h]]}/(E^2+F^2+G^2-a^2+c_1h+\dots+c_lh^l)
$$
is a formal deformation of $C_{\xi}$. If one chooses
$a,c_1,\dots ,c_l$ to be real, then it becomes
the complexification of a 
formal deformation  of the real orbit $C_\xi \cap su(2)$.
\smallskip
To go to the noncompact form it is enough to take the basis $\{\tilde
E=iE, \tilde F=iF, G\}$. The deformed algebra is
$$
U_{[[h]]}/(-\tilde E^2-\tilde F^2+G^2-a^2+c_1h+\dots+c_lh^l).
$$
A basis for $U/I_0$ is
$$
\{g^m\tilde e^n\tilde
f^\mu\}_{m,n=0,1,2\dots\atop \mu=0,1}\;\;.
$$
 The
subalgebra invariant under the automorphism $A$ of Example 2.1, has
instead a basis
$$
\{g^m\tilde e^{2n-m}\tilde f^\mu\}_{m,n=0,1,2\dots\atop \mu=0,1}\;\;.
$$
We can also express this algebra in terms of  the set of commutative generators
$$
v_1=g^2,\quad v_2=\tilde e^2,\quad v_3=g\tilde e,\quad v_4=\tilde f
$$
with relations
$$
v_3^2=v_1v_2,\quad v_1-v_2-v_4^2=a^2.
$$ 
It is clear that this algebra separates the points of the real orbit.
Since the Casimir element is invariant under the automorphism $A$
(extended to $U_h$), it restricts to an automorphism of  $U_h/I_h$.
Analogously to the commutative case, the
subalgebra of $U_h/I_h$ invariant under $A$ can be given in terms of the
generators
$$
V_1=G^2,\quad V_2=\tilde E^2,\quad V_3=G\tilde E ,\quad V_4=\tilde F
$$
and relations
$$
V_3^2=V_1V_2 -h V_3V_4-h^2V_1,\quad V_1-V_2-V_4^2=c(h),
$$
in addition to the commutation relations
$$\eqalign{V_4V_1-V_1V_4 &=h(2V_3) - h^2 V_4,\qquad V_4V_2-V_2V_4=h(2V_3)-h^2V_4,\cr
V_4V_3-V_3V_4&=h(V_1+V_2),\qquad V_3V_1-V_1V_3= -h(2V_1V_4) -h^2V_3,\cr
V_3V_2-V_2V_3&=h(V_4V_2+V_2V_4) +h^2V_3 -h^3V_4,\cr
V_2V_1-V_1V_2&=-h(2V_3V_4)+h^2(V_4^2-V_2-V_1).\cr}
$$
\bigskip
\vfill\eject
{\bf 4. Geometric quantization of $S^2$.}
\bigskip

The subject of geometric quantization is a very vast one and we do not 
intend  to make a review here. Many excellent
reviews exist in the literature (see for example [Pu], [Vo]). We will
try to explain only what is needed to understand the geometric
quantization of our particular case, $S^2$. Some of the results we
exhibit here date back to [So]. We will follow closely the scheme of
[Vo], because there the importance of constructing the algebra of
observables is emphasized.

\medskip

Consider a classical system with phase space $X$ and a group $G$ of symmetries. 
This means that $G$ is a group of
symplectomorphisms of the symplectic manifold $X$,
$$ g\in G,\quad  g:X\mapsto
X \quad\hbox{satisfying}\quad g^*\omega =\omega,
$$
 where $\omega$ is the
symplectic form on $X$.   The
Hamiltonian is a $G$-invariant function, that is, $gH=H$, so $G$ is a
group of symmetries of the equations of motion.

We want to find  a quantization of the classical system that preserves
the symmetry under the group $G$.
The goal of geometric quantization is to construct the Hilbert space $\H$ and
the algebra of quantum observables $\A_h$ acting on $\H$
using only the geometrical elements of
 the classical system. This construction should be ``natural'', that is,
the action of $G$ on $X$ as symplectomorphisms should induce
a unitary representation  of $G$ on $\H$ and an action of $G$ on $\A_h$.
This action should reduce to the conjugation by the unitary representation 
on the operators on $\H$ representing the elements of $\A_h$.
\medskip
{\sl Integral orbit data.}
\bigskip

Let $\xi\in \go^*$ and let $G_\xi$ the isotropy group
of $\xi$ and $\go_\xi$ the corresponding Lie algebra. it is clear that
for $Z\in\go_\xi$, $\hbox{ad}^*_Z\xi=0$, which implies
$$
\xi([Z,Y])=0, \quad \forall \; Y\in \go.\eqno(4.1) 
$$
  Suppose that we have a character $\tau$ of $G_\xi$ satisfying
$$
\tau(\hbox{e}^X)=\hbox{e}^{i\xi(X)},\quad Z\in \go_\xi.
$$
 Such character is called an integral orbit datum. Notice that property
(4.1) is essential. Also, $\xi$ must be such that $\xi(Z)=2\pi m,\; m\in
\Z$ whenever $e^Z=\hbox{Id}$.

>From an integral orbit datum we can construct a unitary representation
of $G$ by induction. We consider the  induced  vector
bundle $E(G/G_\xi, \C_\tau) = (G\times\C)/\tau$, where the equivalence
relation is given by
$$
 (g, v)\quad\approx \quad (gh^{-1},
\tau(h)v), \quad h\in G_\xi.
$$
 We can describe the sections on this
bundle by functions $f:G\mapsto \C$
 satisfying
$$
f(gh)=\tau(h)^{-1}f(g)\eqno(4.2).
$$

By considering the compactly supported sections, and from the fact that
there is a a $G$-invariant measure on $G/G_\xi$ the construction of the
 Hilbert space is straightforward, with bilinear form $$
<f_1,f_2>=\int_{G/G_\xi}f_1\bar f_2.
$$

The problem is that this representation is not necessarily irreducible.
 Nevertheless, in many cases (like for elliptic orbits),
it is possible to restrict naturally the
space of sections  (4.2) to
an irreducible component. We are then interested in computing the integral
orbit data for SU(2).

\medskip

The Lie algebra of SU(2)
 is spanned by the matrices
 $$
 G={i\over
2}\sigma_3,\quad E={i\over 2}\sigma_2,\quad F={i\over 2}\sigma_1
$$
 with
$$
 \sigma_1=\pmatrix{0&1\cr 1&0}, \quad\sigma_2=\pmatrix{0&-i\cr i&0},
\quad \sigma_3=\pmatrix{1&0\cr 0&-1},
 $$
 and commutation relations \footnote*{The spin operators which are used
in physics are given by $G'=-i\hbar G, \; E'=-i\hbar E,\; F'=-i\hbar F$.
 We can reintroduce $\hbar=h/2\pi$ in the analysis
with this rescaling, the multiplication by $-i$ changing a
representation by antihermitian operators of SU(2)  to 
hermitian operators.}
 $$
[E,F]=G,\quad [F,G]=E,\quad [G,E]=F
 $$
 Consider $\xi_a\in \go^*$ such
that $\xi_a(xE+yF+zG)=az$. The isotropy group is
 $$
G_{\xi_a}=\{\hbox{e}^{zG}, z\in \R\}= \{\pmatrix{\hbox{e}^{iz/2}&0\cr
0&\hbox{e}^{-iz/2}}, z\in \R\}
 $$
 with Lie algebra
$\go_{\xi_a}=\hbox{span}\{G\}$. If $z=4\pi n$, $n\in \Z$, then
${e}^{zG}=\hbox{Id}$, so in order to have an integral orbit datum,
 $$
\xi_a(4\pi nG)=4\pi na\in 2\pi\Z \quad \forall\; n,  
 $$
 which is possible if and only if $a\in \Z/2$.
 
 The Cartan-Killing form allows the identification of $\go$ and $\go^*$,
also intertwining the adjoint and coadjoint representations. It is given by
$$
<X,Y>=-{1\over 2}\hbox{Tr}(\hbox{ad}X\hbox{ad}Y),\quad X,Y\in \go
$$
that is,
$$
<E,E>=<F,F>=<G,G>=1
$$
and the rest 0. So $\xi_a\approx aG$, and the orbit is given by the
Casimir polynomial
 $$
C=x^2
+y^2 +z^2=a^2
$$
We conclude that only orbits with half integer radius have
integral orbit data. We will denote by $\tau_m$ the corresponding
integral orbit datum, $\tau_m(e^{zG})=(e^{iz/2})^m$

It is easy to convince oneself that the representation in the space of
functions (4.2) is far too large to be irreducible. To overcome this
problem we need to further restrict  the space of sections. We will do
that with the help of a {\it complex polarization}.
\medskip

{\sl   Complex polarization and Hilbert space.}

\medskip

Elliptic orbits have a $G$-invariant complex structure. We define this
complex structure following [Vo]. From now on we use the identification
between $\go$ and $\go^*$ given by the Cartan-Killing form, so we will
use alternatively $\xi=\xi_X\in \go^*$ with $X\in \go$.
\medskip
{\bf Theorem 4.1}.
{\it Let $X\in \go$ be such that ad$_X$ has only imaginary eigenvalues.
  Let $\gb$ be the complexified Lie algebra of $\go$ and let $\gb^t$ ($t\in
\R$) be the $t$-eigenspace of ad$_{iX}$. Then
$$
 \gb=\sum_{t\in\R}\gb^t, \quad
(\go_X)_c=\gb_X=\gb^0
$$
is a gradation of $\gb$. We define  $$ \p_X
=\sum_{t\geq 0}\gb^t, \quad \n_X =\sum_{t> 0}\gb^t. $$ The following
properties are satisfied

\noindent a. $\gb^s$ and $\gb^t$ are orthogonal unless $s=-t$. 

\noindent b. $\bar{\gb^s}=\gb^{-s}$. (Bar means complex conjugation with
respect to the real form $\go$).

\noindent c. The adjoint action of $G_X$ preserves $\gb^t$.}

\medskip

$\gb/\gb_X\approx T_{\xi_X}(G\cdot\xi_X)_c$ is the complexified  tangent space at the
identity coset. The $G$-invariant complex structure can be characterized by requiring
that $\p_X/\gb_X$ is the antiholomorphic tangent space at the identity coset.

\medskip

Let us write down the standard complex structure on $S^2$ to relate it
with this formalism. Let  $V=xE+yF+zG=x\partial_x+y\partial_y+z\partial_z\in
\go$. We take a representative $aG$ for the orbit of radius $a$,
$$
x^2 +y^2+z^2=a^2.
$$
Stereographic coordinates are given in terms of the embedding coordinates
by
$$
V_1=S^2-\{(0,0,-a)\}, \quad u_1={ax \over z+a},\quad v_1={ay\over z+a}
$$

$$
V_2=S^2-\{(0,0,a)\}, \quad u_2={ax \over z-a},\quad v_2={ay\over z-a}.
$$
The action of SU(2) is the one induced by the adjoint representation of SU(2).

Let $x_1:U_1\longrightarrow \C,\;x_2:U_2\longrightarrow \C$
  be the projective coordinates for the complex projective space
$\P^1=U_1\cup U_2$. If we identify
$$
x_1\equiv -v_1+iu_1, \quad x_2\equiv -v_2-iu_2,
$$ 
we obtain a diffeomorphism
$S^2\approx \P^1$. This gives to $S^2$ the complex structure mentioned
above.  For this particular choice, 
the action of SU(2)   obtained from the three dimensional
representation restricted to $S^2$ coincides with the one obtained
from the fundamental representation with  the projective structure.

\medskip
We write now the complexification of su(2), sl(2,$\C$), in the standard basis
$$
H=-i2G,\quad
X=E-iF,\quad Y=-E-iF.
$$
 The eigenvalues  of $iaG$ are
$ -a, 0, +a$ and the
corresponding eigenspaces are
$$
\gb^0=\hbox{span}\{G\}, \quad \gb^a= \hbox{span}\{Y\},
\quad \gb^{-a}=\hbox{span}\{X\}.
$$
 The tangent
space at the North pole ($x=y=0, z=a$) is spanned by $\partial_x,
\partial_y  \in \go/\go_{aG}$ and in terms of the stereographic
coordinates,
$$
\partial_x={1\over 2}\partial_{u_1},\quad \partial_y={1\over 2}\partial_{v_1}.
$$
In the complexified tangent space,
$$
X=\partial_x-i\partial_y={i\over 2}(-\partial_{v_1}-i\partial_{u_1}),\quad
Y=-\partial_x-i\partial_y={i\over 2}(-\partial_{v_1}+i\partial_{u_1}), 
$$
and since the complex coordinate is $x_1=-v_1+iu_1$,
$$
\gb^a=\hbox{span}\{Y\}=\hbox{span}\{\partial_{\bar x_1}\}.
$$
\medskip
{\bf Definition 4.2}.
 A {\it $G$-invariant  complex polarization} is a lagrangian subspace of the 
complexified tangent bundle  at $\xi$,  $ T_{\xi}(G\cdot\xi)_c \approx
\gb/\gb_\xi$.
\medskip
 We remind that a subspace is a lagrangian subspace if the symplectic form
is 0 on that subspace and its dimension is half the dimension of the
symplectic manifold. Because of property  a in Theorem 4.1, $\p_X/\gb_X$ is a
lagrangian subspace and then a complex polarization.

Consider now an integral orbit datum, $\tau$. One can prove that $d\tau$
extends to a representation $\phi$ of $\p_X$ .This extension
satisfies $\phi|_{\n_X}=0$. The induced bundle associated to the
character $\tau$, $E(G/G_X, \C_\tau)$ has also a complex structure and
the holomorphic sections are characterized by
$$
 Z.f=-\phi(Z)f \qquad Z\in \p_X.\eqno(4.3)
$$
where  $f:G\mapsto \C$ satisfies  $ f(gh)=\tau(h)^{-1}f(g),
\;\; g\in G,\; \; h\in G_X$. We will see that in our case this
constructions gives directly the Hilbert space. For other groups,
further corrections are needed.

It is easy to see that for SU(2) the principal  bundle
$E(\hbox{SU(2)/U(1)}, \hbox{U(1)})$ is only a reduction  of the principal bundle
given by the natural projection
$$
\pi : \C^2-\{0\}\mapsto S^2\approx \P_1
$$
that we call $\Theta( S^2, \C^*)$. The corresponding associated bundles
by the representation $\tau_m$ (extended to $\C^*$), will be denoted by
$E(m)$, $\Theta(m)$. $\Theta(m)$ is an holomorphic vector bundle, whose
sections satisfy (4.3), which in this case is simply
$$
\partial_{\bar x_1}f=0.
$$
 Line bundles over $S^2$ are well studied. A holomorphic section on
 $\Theta(m)$
$$
s:\P_1\mapsto \Theta(m)/\quad \pi\circ s=\hbox{id}_{\P_1},
$$
can be given in terms of a function
$$
\tilde s:\C^2-\{0\}\mapsto \C_m
$$
($(\lambda,\rho) \in \C^2-\{0\}$)  satisfying $\tilde s(\lambda\cdot \gamma, \lambda\cdot \rho)=
\lambda^m\tilde s(\gamma,\rho)$ where $\tilde s$
is a homogeneous polynomial in two variables of degree $m$. The group
SU(2) naturally acts on this space of sections, constituting the
$(m+1)$-dimensional (unitary) irreducible representation of SU(2).

We see that geometric quantization associates  quite naturally to the
orbit a Hilbert space
where the group $G$ acts. The last step now is to find the algebra of quantum observables.

\medskip

{\sl   Quantum observables.}

\medskip

Following [Vo], the algebra of observables is the algebra of ``twisted differential operators''
[Vo]  on sections
of the bundle given by the polarization (real or complex). These operators are
endomorphisms of the space of sections of the bundle satisfying certain conditions
(which make plausible the name of ``differential operators''). We will not give here the
general definition, but we will work  with the SU(2)-bundles using the
description given above.
\bigskip
Consider the space of functions $f:\C^2 -\{0\}\longrightarrow \C$, and
$(\c, \r)$ global coordinates on $\C^2 -\{0\}$. Consider the
algebra of differential operators generated by the elements
$$
\c\partial_{\c},\quad \c\partial_{\r},\quad \r\partial_\c, \quad \r\partial_\r
$$
We denote this algebra by ${\cal D}$. It is a filtered algebra (each of
the elements above has degree 1).

The algebra of twisted differential operators on $\Theta(m)$ is
$$
{\cal D}_m={\cal D}/(D-m\hbox{Id})
$$
where $D=\c\partial_{\c}+\r\partial_\r$ is an element in the center
of  ${\cal D}$.

We want to give a presentation for ${\cal D}_m$ and compare it to the
algebra $U_h/I_h$ obtained in section 3.
\bigskip

Consider now ${\cal U}$ the universal enveloping algebra of the Lie
algebra su(2)$^\C\approx$sl(2,$\C$). Let $\{X,Y,H\}$ be  the
standard basis of sl(2,$\C$) (Example 3.16),

\medskip

{\bf Lemma 4.3}.
{\it The filtered algebra homomorphism $p:{\cal
U}\longrightarrow{\cal D}$, given by
$$
p(X)=-\c\partial_{\r},\quad p(Y)=-\r\partial_\c,
\quad p(H)=-\c\partial_{\c}+\r\partial_\r.
$$
is injective.
}
\smallskip
{\sl Proof.}
Notice that ${\cal D}$ acts on the space $P_m=
\{$homogeneous polynomials of degree $m\}$. We denote by $R_m:{\cal
D}\longrightarrow \hbox{End}(P_m)$  this representation. Notice that
$\tilde R_m=R_m\circ p$ is the $m+1$-dimensional irreducible representation of
su(2). Since we have that $\tilde R_m (Z)=0 \;\;Z\in {\cal U}\;\;\forall\; m\Rightarrow
Z=0$ [HC], it follows  that $p$ is an injective map.

\medskip

{\bf Lemma 4.4}. 
{\it
$$
 {\cal D}\cong {\cal U}\otimes\hbox{span}\{D\}/(C-{D\over 2}({D\over 2}+1))
$$
where $C={1\over 2}(XY+YX+{1\over 2}H^2)$ is the Casimir element in ${\cal U}$.
}
\smallskip
{\sl  Proof.}
Define the Lie algebra homomorphism
$$
{\cal U}\otimes \hbox{span}\{D\} \mapright{S} {\cal D}
$$
as $S(W \otimes D)=p(W)D$. Since  $\{p(X),p(Y),p(H), D\}$ generate ${\cal D}$,
$S$ is surjective. We want to show that ker$S$=$I$, where
$I=(C-{D/ 2}({D/ 2}+1))$.
One can check directly that $I \subset \hbox{ker}S$.
We prove $\hbox{ker}S \subset I$ by contradiction.

Observe first that any element $P \in {\cal
U}\otimes\hbox{span}\{D\}/(C-{D\over 2}({D\over 2}+1))$ can be written
as $AD+B$. In fact, let $P=\sum_{k=0}^NA_kD_k$. By induction on $N$.
The cases of $N=0, 1$ are obvious. Let $N>1$. 
$$
P=A_ND^{N-2}(4C-2D)+\sum_{k=0}^{N-1}A_kD_k
$$
By induction we have our result.

\medskip

Let $P_{N-1}= B_1D+B_0$ be a non zero element in ker$S$ that is not in $I$.
Let us construct the combination
$$
P'_{N-1}= B_1P_1+ {1\over 4}P_{N-1}=({1\over 4}B_0-{1\over 2}B_1)D +B_1C.
$$
it is clear that $P'_{N-1}$ doesn't belong to $I$ unless it is identically
0, that is,  $B_0=B_1=0$. In this case  $P_{N-1}$ is also 0, against the hypothesis.
So $P'_{N-1}$ is in $\hbox{ker}(S)$ and not in $I$.  Let us construct now
the combination
$$
P_N=({1\over 4}B_0-{1\over 2}B_1)P_{N-1}-P'_{N-1}B_1=
{1\over 4}B_0^2 -{1\over 2}B_1B_0 -B_1^2C.
$$
Since $P_N \in $ ker$S$ and $P_N$ does not contain $D$, by the injectivity
of $p$ we must have $P_N=0$, that is
$$
{1\over 4}B_0^2-B_1^2C ={1\over 2}B_1B_0  .
$$
Similarly if we construct	
$$
P'_N=P_{N-1}({1\over 4}B_0-{1\over 2}B_1)-P'_{N-1}B_1=
{1\over 4}B_0^2 -{1\over 2}B_0B_1 -B_1^2C.
$$
$P'_N$ must also be 0, so we have that
$$
{1\over 4}B_0^2-B_1^2C ={1\over 2}B_0B_1  .
$$
It follows that $B_1$ and $B_0$ commute. Lets us rewrite any of these
two relations as
$$
(B_0-B_1)^2=(4C+1)B_1^2.\eqno(4.4)
$$
We show that this relation cannot be satisfied unless $B_0=B_1=0$ and 
this will be a contradiction.
Consider the homomorphism from the (filtered) enveloping algebra to the
(graded) symmetric algebra, given by the natural projections
$$
\pi_n: {\cal U}^{(n)}\longrightarrow S^n={\cal U}^{(n)}/{\cal U}^{(n-1)}
$$
and project (4.4) to the symmetric algebra (isomorphic to the
polynomial algebra). It is obvious that the polynomial $\pi_n(4C+1)$
 is not the square  of another polynomial. It follows that
(4.4) cannot be satisfied unless $B_0=B_1=0$.

{\bf Theorem 4.5}.
$$
{\cal D}_m= {\cal U}/(C-{m\over 2}({m\over
2}+1)\hbox{Id}).
$$
\smallskip

{\sl Proof.} Immediate from the definition of ${\cal D}$ and the lemma (4.2).
\medskip

We now want to make an explicit 
comparison with the result of deformation quantization, let us
make the rescaling
$$\tilde X\mapsto \hbar X,\quad  \tilde Y\mapsto \hbar Y,
\quad
\tilde H\mapsto \hbar H,\quad \tilde D=\hbar D.
\qquad\qquad  (4.4)$$
In what follows, $\hbar$ is a number, not an indeterminate; so we are
comparing the geometric quantization with the specialization for a value
of $\hbar$ of the deformation of the polynomial algebra obtained in section
3. Notice that with this rescaling we obtain  a family of isomorphic Lie algebras
$$
[\tilde H,\tilde X]=\hbar 2\tilde X,\quad [\tilde H,\tilde Y]=-\hbar 2\tilde Y,
 \quad [\tilde X,\tilde Y]=\hbar \tilde H.
$$
(and $\tilde D$ in the center) except for $\hbar\mapsto 0$
(while keeping the generators constant) in which the algebra
 becomes abelian. ${\cal U}_\hbar$ is the enveloping algebra of the Lie algebra
for each value of $\hbar$.

The Casimir operator is
$$
\tilde C={1\over 2}(\tilde X\tilde Y+\tilde Y\tilde X +{1\over 2}\tilde H^2).
$$
Using (4.4), the corresponding ideal in  ${\cal U}_\hbar$ is

$$
(\tilde C-l(l+\hbar )),\qquad l=\hbar {m/2}.
$$
It is enough to take $c(\hbar)=l(l+\hbar )$ to obtain the result of
section 3.

Since $l$ is the eigenvalue of the central element $D/2$ in the
corresponding representation, taking the limit $\hbar\mapsto 0$ and keeping
the generators constant (abelian Lie algebra) is equivalent to take
$m\mapsto \infty$. In the physical picture one says that the
classical limit corresponds to large quantum numbers.

\medskip

We want to make the following observations. By choosing different
polynomials $c(h)$ and different values of $h$ 
we obtain  that the specialized $\C$-algebras in general are
not isomorphic. In fact, it is a known
result (see [Va1]) that ${\cal U}/(C-\mu 1)$ has no finite dimensional
representations when $\mu$ is not rational, 
hence different values of $\mu$ (that is of $c(h)$) may give non
isomorphic algebras.

We also  want to remark that our deformation quantization not only
gives a subalgebra that can be specialized for any value of $h$ (namely
the subalgebra of elements that have coefficients that are polynomials
in $h$), but in the special case of $SU(2)$, $SL(2,\C)$, when $h$
is taking certain values, realizes the subalgebra as a concrete algebra
of differential operators on the space of sections described above.

Finally, comparing with the approach of [BBEW], it is easy to see that
the subalgebra of observables with converging star product is the same
as the one we obtain, that is, the algebra of polynomials on the
algebraic manifold.

\bigskip

\centerline{{\bf Acknoledgements}}
\bigskip
 We wish to thank especially Prof. Varadarajan for his help
 during the preparation of this paper. We also wish to
thank A. Brown,  S. Ferrara, D. Gieseker and  G. Mess for helpful discussions.
 
\bigskip

\centerline{\bf REFERENCES}
\bigskip

\item {[ACG]} D. Arnal, M. Cahen and  S. Gutt. {\it Deformations on
   coadjoint orbits.} J. of Geom. and Phys. {\bf 3} (1986), 32-351.

\item{[ALM]} D. Arnal, J. Ludwig, M. Masmoudi. {\it D\'eformations
covariantes sur les orbites polaris\'ees d'un group de Lie}. J. of Geom.
and Phys. {\bf 14} ,  309-331,  (1994).

\item{[Be]} F. A. Berezin. {\it Quantization}. Math. USSR Izvestija,
{\bf 8} 1109-1165 (1974).

\item{[BFFLS]} F. Bayen, M. Flato, C. Fronsdal, A. Lichnerowicz, D.
Sternheimer. {\it Deformation theory and quantization. I and II}. Ann.
Phys., {\bf 111} 1 61-151 (1978).

\item{[Bo]} A. Borel. {\it Linear Algebraic Groups}. Springer Verlag,
(1991).

\item{[BBEW]}  M. Bordemann, M. Brischle, C. Emmrich, S. Waldmann. 
{\it Subalgebras with converging star products in deformation
quantization}. Preprint q-alg/9512019 (1995).

 \item {[CG]} M. Cahen, S. Gutt . {\it Produits $*$ sur les
		orbites des groupes semi-simples
		    de rang 1}, C.R. Acad. Sc. Paris  296 (1983),
		    s\'erie I, 821-823; {\it An algebraic construction
		of $*$ product on the regular
		    orbits of semisimple Lie groups}. In {\sl Gravitation
and Cosmology}. Monographs and Textbooks in Physical Sciences. A volume
in honor of Ivor Robunian, Bibliopolis. Eds W. Rundler and A. Trautman, (1987);
 {\it Non localit\'e d'une
      d\'eformation symplectique sur la
  sph\`ere $S^2$}.  Bull. Soc. Math. Belg. {\bf 36 B} (1984) 207-21.

\item{[CP]} V. Chari,  A. Pressley. {\it A Guide to Quantum Groups}.
Cambridge University Press, (1994).

\item{[DL]} M. De Wilde, P. B. A. Lecomte. {\it Existence of
star-products and of formal deformations  in Poisson Lie algebras of
arbitrary symplectic manifolds}. Lett. Math. Phys., {\bf 7}, 487-496, (1983).

\item{[EK]} P. Etingof, D. A. Kazhdan. {\it Quantization of Lie
bialgebras, I}. Selecta Math., New Series {\bf 2}, 1, 1-41, (1996).

\item{[Fe]}  B. Fedosov,  {\it A simple geometric construction of
deformation quantization}. J. Diff. Geom., {\bf 40}, 2, 213-238, (1994).

\item{[GV]} J. M.  Garc\'{\i}a Bond\'{\i}a and Joseph C. Varilly.
 {\it From geometric quantization to Moyal quantization}.
J. Math. Phys., 36, 2691-2701, (1995).

\item{[Gu]} S.  Gutt. {\it An explicit $*$-product on the cotangent bundle of a Lie group}.
Lett. Math. Phys. 7 (1983) 249-258.

\item{[HC]} Harish-Chandra {\it Representations of a semisimple Lie
group I, II, II}. Trans. Amer. Math. Soc., {\bf 75}, 185-243, (1953),
{\bf 76}, 26-65, (1954), {\bf 76}, 234-253, (1954).

\item{[Ho1]} J. Hoppe. {\it Quantum theory of a massless relativistic surface and a two dimensional bound state problem}. 
MIT PhD thesis, (1982).

\item{[Ho2]} J. Hoppe. {\it On the deformation of time harmonic flows}.
In {\sl Deformation Theory and Symplectic Geometry}, proceedings of the
Ascona meeting, June 1996. Eds. D. Sternheimer, J. Rawnsley and S. Gutt.
Kluwer Academic Publishers (1997). 

\item{[Ki]} A. Kirillov. {\it Elements of the Theory of Representations}.
Springer Verlag, (1975). 

\item{[Ko]} M. Kontsevich. {\it Deformation Quantization of Poisson
Manifolds}. Preprint \break q-alg/9709040, (1997).

\item{[Ks]} B. Kostant. {\it Lie group representations on polynomial rings}.
Amer. J. of Math., {\bf 86}, 271-309, (1964).

 \item{[Pu]} M. Puta. {\it Hamiltonian Mechanical Systems and Geometric
Quantization}. Kluwer Academic Publishers, (1993).

\item{[RCG]} J.Rawsnsley, M. Cahen and S. Gutt. {\it Quantization of
K\"ahler manifolds I: geometric interpretation of Berezin's
quantization}. JGP {\bf 7} no.1, 45-62 (1990).

\item{[RT]} N. Reshetekhin, L. Takhtajan. {\it Deformation Quantization
of K\"ahler manifolds}. \break Preprint, (1999).

\item{[So]} Sorieau. {\it  Structure des syst\'emes dynamiques}. Dunod,
Paris. (1970).

\item{[Va1]} V. S. Varadarajan.  {\it Lie groups, Lie Algebras and their
Representations}. Springer Verlag, (1984).

\item{[Va2]} V. S. Varadarajan  {\it Harmonic  Analysis on Real Reductive
Groups.} Lecture Notes in Mathematics, no. 576. Springer-Verlag, (1977).

\item{[Va3]} Varadarajan, V. S.
{\it On the ring of invariant polynomials on a semisimple
Lie algebra}. Amer. J. Math, {\bf 90}, (1968).

\item{[Vo]} D. Vogan.  {\it The Orbit Method and Unitary Representations for 
Reductive Lie groups}, in {\sl Algebraic and Analytic Methods in Representation
Theory.} Perspectives in Mathematics, Vol 17. Academic Press. (1996).

\end